\documentstyle[12pt]{amsart}
\newcommand{\bC}{{\Bbb C}}
\newcommand{\bP}{{\Bbb P}}

\newcommand{\cA}{{\cal A}}
\newcommand{\cE}{{\cal E}}
\newcommand{\cJ}{{\cal J}}
\newcommand{\cL}{{\cal L}}
\newcommand{\cO}{{\cal O}}

\newcommand{\cQ}{{\cal Q}}
\newcommand{\cR}{{\cal R}}

\newcommand{\tC}{\widetilde{C}}
\newcommand{\tP}{\widetilde{P}}

\newcommand{\trho}{\widetilde{\rho}}
\newcommand{\tSi}{\widetilde{\Sigma}}

\newcommand{\ra}{\rightarrow}
\newcommand{\lra}{\longrightarrow}

\newcommand{\surj}{\hspace{3pt}\to \hspace{-19pt}{\rightarrow} \:\:}
\newcommand{\G}{\Gamma}
\newcommand{\Si}{\Sigma}
\newcommand{\D}{\Delta}
\newcommand{\T}{\Theta}
\newcommand{\tT}{\widetilde{\Theta}}
\newcommand{\half}{\frac{1}{2}}

\newtheorem{proposition}{Proposition}[section]
\newtheorem{lemma}[proposition]{Lemma}

\newtheorem{theoremi}{Theorem}
\newtheorem{conjecture}[theoremi]{Conjecture}

\newtheorem{corollary}[proposition]{Corollary}
\newtheorem{corollaryi}[theoremi]{Corollary}
\newtheorem{propositioni}[theoremi]{Proposition}
\newtheorem{remark}[proposition]{Remark}

\begin{document}

\title{Second order theta divisors on Pryms}
\author{E. Izadi}
\address{Department of Mathematics, Boyd
Graduate Studies Research Center, University of Georgia, Athens, GA
30602-7403, USA}
\email{izadi@@math.uga.edu}

\maketitle

The Schottky problem is the problem of finding necessary and
sufficient conditions for a principally polarized abelian variety
({\em ppav}) to be a product of jacobians of smooth curves.

Let $(P, \Xi)$ be a ppav of dimension $p\geq 4$ with $\Xi$ a symmetric
theta divisor on $P$. The elements of $\G = H^0(P, 2\Xi)$ are
symmetric, hence their multiplicities at the origin are always
even. Let $\G_0
\subset \G$ be the subvector space of $\G$ of sections with
multiplicity at least $2$ at the origin and let $\G_{00} \subset \G$
be the subvector space of sections with multiplicity at least $4$ at
the origin. Also let $|2 \Xi |_{00} \subset |2 \Xi |_0 \subset |2 \Xi
|$ be the linear systems of divisors of zeros of elements of $\G_{00}
\subset \G_0\subset \G$ respectively. It
is well-known that the dimensions of $\G, \G_0$ and $\G_{00}$ are
respectively $2^p$, $2^p -1$ and $2^p -1- \frac{p(p+1)}{2}$ (for this
last we need $(P, \Xi)$ to be indecomposable, see \cite{igusa72} page
$188$ Lemma $11$ and \cite{geemengeer86} Proposition 1.1 page 618).

We have the linear map
\[ \tau_4 : \G_{00} \lra H^0(\bP T_0 P, \cO_{\bP T_0P}(4)) \]
defined by sending a section $s \in \G_{00}$ to the quartic term of
its Taylor expansion at the origin. Let $\cQ_{00}$ denote the linear
subsystem of $|\cO_{\bP T_0P}(4)|$ which is the projectivization of
the image of $\tau_4$. Let $V_{00}$ and $V_{inf,00}$ denote the base
loci of $|2 \Xi |_{00}$ and $\cQ_{00}$ respectively. In
\cite{geemengeer86}, van Geemen and van der Geer proposed a
characterization of the locus of jacobians, made more precise by
Donagi (\cite{donagi88} page 110), in the following form:

\begin{conjecture}
\label{conjsch}
\begin{enumerate}
\item \label{conjjac} If $(P, \Xi)= (JC, \T)$ is the jacobian of a
smooth curve $C$ of genus $g$, then $V_{00}$ is set-theoretically
equal to the reduced surface $C-C := \{\cO_C(s-t) : s,t\in C \}$.

\item \label{conjV00} If $(P, \Xi)$ is not in the closure $\cJ_g$ of
the locus of jacobians in the moduli space of ppav of dimension $g$, then
$V_{00}=\{ 0 \}$ set-theoretically.

\end{enumerate}
\end{conjecture}

Beauville and Debarre proposed an infinitesimal version of
Conjecture \ref{conjsch} (\cite{beauvilledebarre89}):

\begin{conjecture}
\label{conjschinf}
\begin{enumerate}
\item \label{conjjacinf} If $(P, \Xi)= (JC, \T)$ is the jacobian of a
smooth curve $C$ of genus $g$, then $V_{inf,00}$ is,
set-theoretically, the canonical image $\kappa C$ of $C$ in $|\omega_C
|^* = \bP T_0 JC$ where $\omega_C$ is the dualizing sheaf of $C$ (note
that $\kappa C$ is the projectivized tangent cone to $C-C$ at $0$).

\item \label{conjV00inf} If $(P, \Xi)$ is not in $\cJ_g$, then $V_{inf,00}$
is empty.
\end{enumerate}
\end{conjecture}

The first parts of Conjectures \ref{conjsch} and \ref{conjschinf} have
been proved (each with one well-determined exception) by Welters
\cite{welters86}, the author \cite{I1} and Beauville and Debarre
\cite{beauvilledebarre89}.
In \cite{I1}, we also determined the scheme-structures of the base
loci for jacobians. In \cite{I3}, Theorem 4 page 95, we proved the
second parts of Conjectures \ref{conjsch} and \ref{conjschinf} in the
case $p=4$. Beauville, Debarre, Donagi and van der Geer
proved part $2$ of Conjecture \ref{conjsch} for intermediate jacobians
of cubic threefolds, the Prym varieties of ``even'' \'etale double
covers of smooth plane curves (see \cite{BDDG}), and Beauville and
Debarre proved parts 2 of Conjectures \ref{conjsch} and
\ref{conjschinf} for certain ppav isogenous to a product of $p$
elliptic curves (see \cite{beauvilledebarre89} pages 35-38). By
semi-continuity, one deduces from this last result that for a {\em
general} ppav $V_{00}$ is finite and $V_{inf,00}$ is empty.

Let $\cA_p$ be the moduli space of ppav of dimension $p$. In
\cite{I3}, we proved the second parts of Conjectures \ref{conjsch} and
\ref{conjschinf} for $p=4$ by using the fact (proved
in \cite{I3}, Theorem 3.3 page 111) that an element of $\cA_4\setminus
\cJ_4$ is always the Prym variety of an \'etale double cover of a {\em
smooth} curve of genus $5$. From now on we will suppose that $(P,
\Xi)$ is the Prym variety of an \'etale double cover of smooth curves
$\pi : \tC \rightarrow C$ with $C$ non-hyperelliptic of genus $g=p+1$.
Then $(P,\Xi)$ is indecomposable, i.e., it is not the product of two
positive-dimensional ppav (see \cite{mumford74} page 344,
Theorem(d)). There is a natural analogue of the surface $C-C$ for a
Prym variety, namely, the reduced surface
\[\Si := \Si(\pi : \tC \rightarrow C) := \{\cO_{\tC}(s+t-\sigma
s - \sigma t) : s,t \in \tC \} \subset P\subset J\tC,
\]
where $\sigma$ is the involution of the cover $\pi : \tC \rightarrow
C$. Let $\epsilon :\tP\ra P$ be the blow up of $P$ at $0$ with
exceptional divisor $\cE$ and let $\tSi$ be the proper transform of
$\Si$ in $\tP$. Let $L$ be the linear system $|\epsilon^* (2 \Xi) -
4\cE |$ on $\tP$. When $g=5$ there is an involution $\lambda$ acting
on the moduli space of admissible double covers of stable curves of
genus $5$ such that a double cover $\tC\ra C$ and $\lambda (\tC\ra C)$
have the same Prym variety (see \cite{donagi92} page 100 and \cite{I3}
pages 119 and 126). Furthermore, for any fixed
$(P,\Xi)\in\cA_4\setminus\cJ_4$, there is an \'etale double cover
$\tC\ra C$ of a smooth curve $C$ such that $(\tC_{\lambda }\ra
C_{\lambda }) :=\lambda (\tC\ra C)$ is also an \'etale double cover of
a smooth curve $C_{\lambda }$ (see \cite{I3} page 136). In such a
case, put $\Si_{\lambda } :=\Si(\tC_{\lambda }\ra C_{\lambda })$ and
let $\tSi_{\lambda }$ be the proper transform of $\Si_{\lambda }$ in
$\tP$. With these hypotheses, we proved in \cite{I3} that (recall $g-1=p=4$)

\begin{itemize}

\item there is exactly a pencil of elements of $L$ containing
$\tSi$ (see \cite{I3}, 5.7 page 134 and 6.23 page 148),

\item the base locus of this pencil is equal to
$\tSi\cup\tSi_{\lambda }$ as a set, as a scheme if $(\tC\ra C)\neq
(\tC_{\lambda}\ra C_{\lambda})$ (follows from \cite{I3}, 5.7 page
134),

\item the base locus of the restriction $L |_{\tSi}$ is
empty (\cite{I3} pages 139 and 146-147).

\end{itemize}
We generalize this third result to higher dimensional Prym varieties and
give bounds for the dimension of the sublinear system of $L$ consisting of
elements containing $\tSi$:
\vskip15pt

In Section \ref{secproof} we prove
\begin{theoremi}\label{maintheorem}
If $C$ is neither trigonal nor hyperelliptic of genus $\geq 5$, then the
base locus of $L |_{\tSi}$ is empty.
\end{theoremi}
\noindent
For $C$ trigonal, the Prym $(P,\Xi)$ is a jacobian and the support of
the base locus of $L |_{\tSi}$ is easily determined (see Proposition
\ref{rhobase} and Remark \ref{rembase} below).

To prove the Theorem, we use divisors in the linear system $|2 \Xi |_{00}$
which are obtained as intersections with $P \subset J \tC$ of
translates of the theta divisor $\tT$ of $J \tC$ (see Section \ref{newdiv}
below). To our knowledge such divisors have not been used before in the
litterature.

The Prym-canonical curve $\chi C$ is the image of $C$ in
$|\omega_C\otimes\alpha |^*$ by the natural morphism $C\ra
|\omega_C\otimes\alpha |^*$ where $\alpha$ is the square-trivial
invertible sheaf associated to the double cover $\pi :\tC\ra C$. Under
the natural isomorphism $\bP T_0 P\cong |\omega_C\otimes\alpha |^*$,
the curve $\chi C$ is the tangent cone to $\Si$ at $0$. Therefore the
Theorem implies
\begin{corollaryi} \label{maincor} Suppose that $C$ is neither
trigonal nor hyperelliptic, then
\begin{enumerate}
\item the only base point of $|2 \Xi |_{00}$ on $\Si$ is $0$,
\item the linear system $\cQ_{00}$ has no base points on $\chi
C$.
\end{enumerate}
\end{corollaryi}

Our approach leads us to consider the vector space $\G_{00}' :=\{
s\in\G_{00} : s |_{\Si} = 0\}$ with projectivization $| 2\Xi |_{00}' := \{
D \in | 2\Xi |_{00} : D \supset \Si \}$. Then $|2\Xi |'_{00}$ can be
identified with the sublinear system of elements of $L$ containing
$\tSi$. Let $\cQ_{00}'$ be the linear subsystem of $\cQ_{00}$ consisting of
quartic tangent cones at $0$ to elements of $| 2 \Xi |_{00}'$. On a Prym
variety, second order theta divisors which contain $\Si$ can be thought of
as natural generalizations of $2\T$-divisors on JC containing $C-C$. It is
well-known (see \cite{fay84} Theorem 2.5 page 120 and \cite{geemengeer86}
page 265 or \cite{welters86} Proposition 4.8 page 18) that, on a jacobian
$(JC,\T)$, a $2\Theta$-divisor contains $C-C$ if and only if it has
multiplicity at least $4$ at the origin. This does not generalize to Prym
varities: Corollary \ref{maincor} implies that $| 2\Xi |_{00}'$ and
$\cQ_{00}'$ are proper linear subsystems of $| 2\Xi |_{00}$ and $\cQ_{00}$
respectively. More precisely, we prove (see Section \ref{secproof}):

\begin{propositioni}\label{mainprop2}
Suppose $C$ non-hyperelliptic. The dimension of $\G_{00}'$ is at most
$2^p - 4 -\half p (p+1)$ and at least $2^p - 2 - p (p-1)$. The
codimension of $\cQ_{00}'$ in $\cQ_{00}$ is at least $2$.
\end{propositioni}
Note that when $p=4$, we have another proof of the result of \cite{I3}
saying that the dimension of $|2\Xi |_{00}'$ is $1$. We pose

\begin{conjecture}\label{newconj}
Suppose that $(P,\Xi)$ is not a jacobian. Then
\begin{enumerate}
\item the base locus $V_{00}'$ of $| 2 \Xi |_{00}'$ is $\Si$ as a set
if $p \geq 6$,
\item the base locus $V_{inf,00}'$ of $\cQ_{00}'$ is $\chi C$ as a
set if $p \geq 6$.
\end{enumerate}
\end{conjecture}

By Corollary \ref{maincor}, this conjecture implies the second parts of
Conjectures \ref{conjsch} and \ref{conjschinf} for Prym varieties. We have
the following evidence for this conjecture:

Results of Welters and Debarre (\cite{welters872} and \cite{debarre89}, see
Section \ref{prelim} below) imply
\begin{propositioni}\label{mainprop}
The base locus $V_{00}'$ is the set-theoretical union of $\Si$ and,
possibly, some curves and points for a general Prym variety of
dimension $p \geq 16$.
\end{propositioni}

Combined with Corollary \ref{maincor} this implies:

\begin{corollaryi}
If $(P,\Xi)$ is a general Prym variety of dimension $\geq 16$, then
$V_{00}$ has dimension $\leq 1$.
\end{corollaryi}

Results of Debarre imply (see Section \ref{prelim})

\begin{propositioni}\label{propinf00}
For a general Prym variety of dimension $p\geq 8$, the base locus
$V_{inf,00}'$ is set-theoretically equal to $\chi C$.
\end{propositioni}
Combined with Corollary \ref{maincor} this implies (see Section
\ref{prelim} below)
\begin{corollaryi}\label{maincorol}
Part 2 of Conjecture \ref{conjschinf} is true for general Prym
varieties of dimension $p=g-1 \geq 8$.
\end{corollaryi}

We explain why we make Conjecture \ref{newconj} only for $p\geq 6$ and
only set-theoretically. Let $\cR_g$ be the space parametrizing \'etale
double covers $\tC\rightarrow C$ where $C$ is a smooth
non-hyperelliptic curve of genus $g=p+1$. The Prym map is the morphism
$\cR_g\ra\cA_p$ which to a double cover $\tC\ra C$ associates its Prym
variety. Recall (see above) that in case $p=4$, we proved in \cite{I3}
(5.7 page 134) that $V_{00}' =\Si\cup
\Si_{\lambda}$. Since the Prym map is generically injective for $g
\geq 7$ (see \cite{friedmansmith82}), we expect that $V_{00}' =\Si$ as
sets. Now an argument analogous to \cite{I1} (2.9) page 196 shows that, if
$V_{00}' =\Si$ as sets, then $V_{00}'$ is not reduced at $0$ and hence is
not equal to $\Si$ as a scheme.  If the morphism $\rho :\tC^{ (2) }\ra\Si$,
where $\tC^{ (2) }$ is the second symmetric power of $\tC$, is birational,
then, by analogy with the case of jacobians (see \cite{I1}), we can expect
$V_{00}'$ to be reduced at the generic point of $\Si$. If, on the other
hand, the morphism $\rho$ is not birational, then, by a standard
semi-continuity argument, the scheme $V_{00}'$ is not reduced anywhere on
$\Si$. Note that, using a refinement of a theorem of Martens by Mumford
(\cite{ACGH} Theorem 5.2 page 193), one can easily see that if $\tC$ is
neither bielliptic, trigonal, nor a smooth plane quintic, then $\rho$ is
birational. Similarly, we can expect $V_{inf,00}'$ to be equal to $\chi C$
as schemes if $C\cong\chi C$ but $V_{inf,00}'$ will not be equal to $\chi
C$ as a scheme if the morphism $C\ra\chi C$ is not (at least) birational.

On the other hand, Donagi and Smith proved (see \cite{donagismith81})
that the generic fibers of the Prym map have cardinality $27$ when
$g=6$ and Donagi proved (see \cite{donagi92} Theorem 4.1 page 90) that
the Galois group of the Prym map is isomorphic to the Galois group of
the $27$ lines on a cubic surface. So there could be nontrivial
automorphisms acting in the fibers of the Prym map. If there is an
automorphism $\mu :\cR_6\ra\cR_6$ acting in the fibers of the Prym
map, then one would expect $V_{00}'$ to be the union of $\Si =\Si
(\tC\ra C)$ and the surfaces $\Si (\mu (\tC\ra C) ),\Si (\mu^2 (\tC\ra
C) ),$ etc.

\vskip20pt

\begin{center}
{\sc Notation and conventions}
\end{center}

All varieties and schemes are over the field $\bC$ of complex numbers.

For any section $s$ of an invertible sheaf $\cL$ on a variety $X$,
denote by $Z(s)$ the divisor of zeros of $s$. Let $h^i(X,\cL)$
denote the dimension of $H^i( X,\cL)$.

For a divisor $D$ on $C$ (resp. $\tC$), denote by $\langle D\rangle$
the span of $D$ in the canonical space of $C$ (resp. $\tC$).

\section{Preliminaries}\label{prelim}

Let $\pi : \tC \ra C$ be an \'etale double cover of a smooth
non-hyperelliptic curve $C$ of genus $g$. Let $\alpha$ be the point of
order $2$ in $Pic^0 C$ associated to the double cover $\pi$ so that we
have $\pi_*\cO_{\tC }\cong\cO_C\oplus\alpha$. Choose an element
$\beta$ of $Pic^0 C$ such that $\beta^{\otimes 2} = \alpha$ and a
theta-characteristic $\kappa$ on $C$ such that $h^0( C,\kappa )$ and
$h^0(\tC ,\pi^* (\kappa\otimes\beta ))$ are even. Symmetric principal
polarizations on $J\tC = Pic^0\tC$ and $JC = Pic^0 C$ can be defined
as the reduced divisors
\[
\tT := \{ D \in Pic^0 \tC : h^0(\tC, D \otimes \pi^* (\kappa \otimes
\beta) ) > 0 \}
\]
\[
\T := \{ D \in Pic^0 C : h^0(C, D \otimes \kappa ) > 0 \} \: .
\]
With these definitions, the inverse image of $\tT$ by the morphism
$\pi^* : JC\ra J\tC$ is the divisor $\T_{\beta} +\T_{\beta^{-1}}$
where $\T_{\beta}$ is the translate of $\T$ by $\beta$. The Prym
variety $( P,\Xi )$ of the double cover $\pi :\tC\ra C$ is defined by
the reduced varieties
\[ P := \{ E \in Pic^0 \tC : Nm(E) \cong \cO_C, h^0(\tC, E \otimes
\pi^*(\kappa \otimes \beta))
\equiv 0 \hbox{ mod }2 \} \: , \]
\[ \Xi := \{ E \in P : h^0(\tC, E \otimes \pi^* (\kappa \otimes \beta))
> 0 \}
\]
where $Nm : Pic \tC \ra Pic C$ is the Norm map (see \cite{mumford74}
pages 331-333 and 340-342). As divisors we have $ 2\Xi = P.\tT$. Also
define
\[
P' := \{ E \in Pic^{2g-2} \tC : Nm(E) \cong \omega_C, h^0( \tC, E)
\equiv 0 \hbox{ mod }2 \}
\]
and
\[
\Xi' :=\{ E \in P' : h^0( \tC, E ) > 0\}\; .
\]
For any $E \in P'$, put
\[ \Xi_E := \{ D \in P : D \otimes E \in \Xi' \} \: .
\]
Since $Nm(E)\cong\omega_C$, we have $\omega_{\tC} \otimes E^{-1}
\cong\sigma^* E$. By the theorem of the square, the divisor $\Xi_E +
\Xi_{\sigma^* E}$ is in the linear system $|2 \Xi |$ and, by Wirtinger
Duality (see \cite{mumford74} pages 335-336), such divisors
span $| 2 \Xi |$. Furthermore, if $E \in \Xi'$, then $\Xi_E +
\Xi_{\sigma^* E}$ is in $| 2 \Xi |_0$ and such
divisors span $| 2 \Xi |_0$ (by Wirtinger duality, the span of such
divisors is the span of $\phi(\Xi)$ where $\phi : P\ra |2\Xi |^*$ is
the natural morphism, this span is a hyperplane in $|2\Xi |^*$ which can
therefore be identified with $| 2\Xi |_0$ by Wirtinger duality, also
see \cite{welters86} page $18$).

We now explain how Proposition \ref{mainprop} follows from results of
Welters and Debarre and how Proposition \ref{propinf00} follows from
results of Debarre. In this paragraph only, suppose that $(P,\Xi)$ is
a {\em general} Prym variety. Then $E$ is an element of the singular
locus $Sing(\Xi')$ of $\Xi'$ if and only if $h^0(\tC, E)\geq 4$
(\cite{welters872} page 168). Therefore, for every $E\in Sing(\Xi')$
and $(p,q )\in\tC^2$, we have $h^0 (\tC, E\otimes\cO_{\tC} (p +q
-\sigma p -\sigma q)) > 0$. Hence $\Si\subset\Xi_E$ and $\Xi_E
+\Xi_{\sigma^* E}$ is in $| 2\Xi |_{00}'$. Since $\Xi$ is symmetric,
the tangent cones at $0$ to $\Xi_E$ and $\Xi_{\sigma^* E}$ are
equal. It follows from these facts and the irreducibility of
$Sing(\Xi')$ for $p\geq 6$ (see \cite{debarre89} Th\'eor\`eme 1.1 page
114) that $V_{00}'$ is set-theoretically contained in $\cap_{E\in
Sing(\Xi')} \Xi_E$ and $V_{inf,00}'$ is set-theoretically contained in
the intersection of the tangent cones at $0$ to $\Xi_E$ for $E$ a
point of multiplicity $2$ on $\Xi'$. This latter can be rephrased as:
$V_{inf,00}'$ is set-theoretically contained in the intersection of
the tangent cones to $\Xi'$ at its points of multiplicity $2$. It is
proved in \cite{welters872} (Theorem 2.6 page 169) that, for $p\geq
16$, the intersection $\cap_{E\in Sing(\Xi')} \Xi_E$ is the
set-theoretical union of $\Si$ and, possibly, some curves and
points. This proves Proposition \ref{mainprop}. By \cite{debarre89}
Th\'eor\`eme 1.1 page 114, the tangent cones to $\Xi'$ at its points
of multiplicity $2$ generate the space of quadrics containing $\chi C$
for $p\geq 6$ and, for $p\geq 8$ (see \cite{debarre89} Corollaire 2.3
page 129, also see \cite{lazarsfeld95} and
\cite{langesernesi96}), the Prym-canonical curve $\chi C$ is cut out
by quadrics. This proves Proposition \ref{propinf00}.

\section{Pull-backs of divisors to $\tC^2$}

Let $\rho : \tC^2 \ra P\subset J\tC$ be the morphism
\[ \rho : (s,t) \longmapsto [s,t] := \cO_{\tC}(s+t-\sigma s-\sigma t) \]
so that $\Si$ is the image of $\tC^{2}$ by $\rho$. The morphism $\rho$
lifts to a finite morphism $\trho :\tC^2\surj\tSi\subset\tP$. We have

\begin{lemma}
For any $E \in \Xi'$ such that $h^0(\tC,E) = 2$ and $|E|$ contains reduced
divisors, the inverse image of $\Xi_E$ in $\tC^2$ is the divisor
\[
\rho^*\Xi_E = D_{\sigma^* E} +\Delta'
\]
where $D_E$ is the reduced divisor
\[
D_E := \{ (p,q) : h^0( \tC, E\otimes\cO_{\tC}(-p-q) ) > 0 \}
\]
and $\Delta'$ is the ``pseudo-diagonal'' of $\tC^2$, i.e., the reduced
curve $\Delta' := \{ (p, \sigma p) \in \tC^2 \} = \rho^{-1}( 0 )$.
\label{XiEC2} \end{lemma}

{\em Proof :} We have the equality of sets
\[
\rho^*\Xi_E = \{ (p,q) : E \otimes [p,q] \in \Xi' \} 
= \{ (p,q) : h^0(\tC, E \otimes [p,q]) > 0 \} \:\: .
\]

First suppose that $| E |$ has no base points. By \cite{mumford71}
(page 187, Step II and its proof), for any $p \in \tC$, we have
$h^0(\tC, E \otimes \cO_{\tC}(p-\sigma p)) = h^0(\tC, E) - 1$ if
$\sigma p$ is not a base point of $|E|$ and $h^0(\tC, E \otimes
\cO_{\tC}(p-\sigma p)) = h^0(\tC, E) + 1$ if $\sigma p$ is a base point of
$|E|$. Since $|E|$ has no base points and $h^0(\tC, E) = 2$, we have
$h^0(\tC, E \otimes \cO_{\tC}(p-\sigma p)) = 1$ for all $p\in\tC$.
Similarly, $h^0 (\tC, E \otimes [p,q]) = h^0 (\tC,
E\otimes\cO_{\tC}(p-\sigma p)) - 1 =0$ if $\sigma q$ is not a base point of
$| E\otimes\cO_{\tC}(p-\sigma p)|$ and $h^0(\tC, E \otimes [p,q]) =
h^0(\tC, E\otimes\cO_{\tC }( p -\sigma p)) + 1 = 2$ if $\sigma q$ is a base
point of $| E\otimes\cO_{\tC } (p -\sigma p)|$. Let $D_{\sigma p}$ be the
unique divisor of $|E|$ containing $\sigma p$, then $D_{\sigma p} +p
-\sigma p$ is the unique element of $| E \otimes\cO_{\tC} (p-\sigma
p)|$. So, to obtain $h^0(\tC, E \otimes [p,q]) > 0$, we need $\sigma q$ to
be a base point of $|E \otimes\cO_{\tC} (p-\sigma p)|$, i.e., the point
$\sigma q$ appears in $D_{\sigma p} +p - \sigma p$. Therefore, either
$\sigma q = p$ or $h^0(\tC, E \otimes \cO_{\tC}(-\sigma p - \sigma q)) >
0$. So the support of $\rho^* \Xi_E$ is $D_{\sigma^* E} \cup \Delta'$.

An easy degree computation on the fibers of $\tC^2$ over $\tC$ by the
two projections shows then that $\rho^* \Xi_E = D_{\sigma^* E}
+ \Delta'$ as divisors.

If $| E |$ has base points, then its base points are distinct because
$| E |$ contains reduced divisors. Then, by \cite{griffithsharris78}
page 287, the cohomology class of $D_{\sigma^* E}$ is the same as in
the case where $|E|$ has no base points. Hence, since $\rho^* \Xi_E =
D_{\sigma^* E} + \Delta'$ as divisors in the case where $|E|$ has no
base points, we see that $\rho^*\Xi_E - \D'$ and $D_{\sigma^* E}$ have
the same cohomology class. Since the support of $D_{\sigma^* E}$ is
clearly contained in the support of $\rho^*\Xi_E -\D'$ and these are
both effective divisors, they are equal.
\qed

\begin{remark}
Suppose that $E$ as above can be written as $E = \pi^* M\otimes
\cO_{\tC}(B)$ where $B$ is the base divisor of $E$ and $M$ is an invertible
sheaf on $C$ with $h^0( C, M) = 2$. Then $D_{\sigma^* E} =\D' + D'$ for
some effective divisor $D'$. Therefore $\rho^*\Xi_E = D' + 2\D'$ which
agrees with the fact that in such a case $E\in Sing (\Xi')$ so that $\Xi_E$
is singular at $0$ (see \cite{mumford74} pages 342-343).
\end{remark}

Let $\Delta \subset \tC^2$ be the diagonal, let $p_i :\tC^2\ra\tC$ be
the projection onto the $i$-th factor and let $\omega_{\tC^2 }$ be the
canonical sheaf of $\tC^2$. Then $\omega_{\tC^2 }\cong
p_1^*\omega_{\tC }\otimes p_2^*\omega_{\tC }$ and we have K\"unneth's
isomorphism $H^0 (\tC^2 ,\omega_{\tC^2 } )\cong H^0 (\tC ,\omega_{\tC
} )^{\otimes 2}$. Let $I_2(\tC) \subset S^2H^0(\tC, \omega_{\tC}
)\subset H^0(\tC ,\omega_{\tC})^{\otimes 2} = H^0(\tC^2,\omega_{\tC^2
})$ be the vector space of quadratic forms vanishing on the canonical
image $\kappa \tC$ of $\tC$. Fix an embedding $H^0(\tC^2
,\omega_{\tC^2 }\otimes\cO_{\tC^2 }( - 2 \Delta) )\subset
H^0(\tC^2,\omega_{\tC^2 } )$ obtained by multiplication by a nonzero
global section of $\cO_{\tC^2 }( 2 \Delta)$ (note that
$h^0(\tC^2,\cO_{\tC^2 }( 2 \Delta)) = 1$ because $\D$ has negative
self-intersection, hence any two such embeddings differ by multiplication
by a constant). Then it is easily seen that
\[
I_2(\tC) = S^2H^0(\tC, \omega_{\tC}) \cap H^0(\tC^2
,\omega_{\tC^2 }\otimes\cO_{\tC^2 }( - 2 \Delta) )\subset
H^0(\tC^2,\omega_{\tC^2 } ) \: .
\]
Similarly fix embeddings $H^0(\tC^2 ,\omega_{\tC^2
}\otimes\cO_{\tC^2 }( - 2 \Delta - 2\D') )\subset H^0(\tC^2 ,\omega_{\tC^2
}\otimes\cO_{\tC^2 }( - 2 \Delta) )\subset H^0(\tC^2,\omega_{\tC^2
}\otimes\cO_{\tC^2} (- 2 \Delta+2 \Delta' ) )$. For $E\in\Xi'$ such
that $h^0(\tC, E)=2$, the linear system $|E|$ has no base points and
$E\not\cong\sigma^* E$, it is well-known (see, e.g.,
\cite{ACGH} page 261) that $q_E :=\cup_{D
\in |E|} \langle D \rangle =\cup_{D \in |\sigma^*E|} \langle D\rangle$
is a quadric of rank $4$ whose rulings cut the divisors of $|E|$ and
$|\sigma^*E|$ on $\tC$. We need the following

\begin{lemma}\label{Zsqs}
\begin{enumerate}
\item\label{restrictC2} We have
\[
\rho^* \cO_P (2 \Xi) \cong \omega_{\tC^2 }\otimes\cO_{\tC^2}
(- 2 \Delta+2 \Delta' ),
\]
\[
\rho^*\G_0\subset I_2 (\tC )\subset H^0( \tC^2
,\omega_{\tC^2 }\otimes\cO_{\tC^2}( - 2\Delta ) )\subset H^0( \tC^2,
\omega_{\tC^2 }\otimes\cO_{\tC^2}(- 2
\Delta+2 \Delta' ) )\; ,
\]
and
\[
\rho^*\G_{ 00 }\subset I_2 (\tC )\cap H^0( \tC^2,\omega_{\tC^2
}\otimes\cO_{\tC^2} ( - 2\Delta - 2 \Delta' ) )\subset H^0( \tC^2
,\omega_{\tC^2 }\otimes\cO_{\tC^2} ( - 2\Delta + 2\Delta' ) ).
\]

\item\label{Zs} For any $s \in \Gamma_0$, let $q(s)$ be the quadric in
$|\omega_{\tC }|^*$ with equation $\rho^* s\in I_2 (\tC
)$. Then, in $\tC^2$, the zero locus of $\rho^* s\in H^0 (\tC^2
,\omega_{\tC^2 }\otimes\cO_{\tC^2}( - 2\Delta +2 \Delta' ) )$ is
\[
Z(\rho^* s) = Z_{q(s)} + 2 \Delta'
\]
where $Z_{q(s)}$ is a divisor with support
\[
\{ (p,q) \in \tC^2 : \langle p+q \rangle \subset q(s) \}.
\]
If $s\in\G_0$ is general, then $Z_{q(s)}$ is reduced. If $Z(s) = \Xi_E
+ \Xi_{\sigma^* E}$ for some $E\in\Xi'$ such that $h^0(\tC, E) = 2$,
then $q(s) = q_E :=\cup_{D\in | E |}\langle D\rangle$.

\item \label{tangcone} If $s \in \Gamma_0 \setminus \Gamma_{00}$, then
\[
q(s) \cap (\bP T_0 P = | \omega_C \otimes \alpha |^*) \subset\bP T_0
J\tC = |\omega_{\tC} |^*
\]
is the projectivized tangent cone $\tau_{Z(s)}$ to $Z(s)\subset P$ at $0$.

\item\label{multDs}
For any $s\in\G_0$, the multiplicity of $\rho^*s$ at the generic point
of $\Delta'$ is even $\geq 2$ and if $\rho^*s$ vanishes on $\Delta'$
with multiplicity $\geq 4$, then either $s
\in \Gamma_{00}$ or $\tau_{Z(s)}$ contains the Prym-canonical curve $\chi
C$.
\end{enumerate}

\end{lemma}

{\em Proof :  }\ref{restrictC2}. Let $E$ be an invertible sheaf of
degree $2g-2$ on $\tC$ such that $h^0(\tC, E)=2$ and $|E|$ has no base
points. Let $s_1$ and $s_2$ be two general sections of $E$. Then
\[ s_1 \otimes s_2 - s_2 \otimes s_1
\in \Lambda^2 H^0(\tC, E) \subset H^0(\tC,E)^{\otimes 2}
\]
and, as in the the proof of Lemma \ref{XiEC2}, it is easily seen that $Z(s_1
\otimes s_2 - s_2 \otimes s_1) = D_E +\Delta$.

Now, supposing that $E \otimes \sigma^*E \cong \omega_{\tC}$, from the natural
map
\[
\psi_E : H^0(\tC, E) \otimes H^0(\tC, \sigma^*E) \lra H^0(\tC, \omega_{\tC})
\]
we obtain the map
\[
\begin{array}{rcl}
H^0(\tC, E)^{\otimes 2} \otimes H^0(\tC, \sigma^*E)^{\otimes 2} & \lra &
H^0(\tC, \omega_{\tC})^{\otimes 2}\cong H^0(\tC^2, \omega_{\tC^2}) \\
t_1\otimes t_2\otimes\sigma^* u_1\otimes\sigma^* u_2 & \longmapsto &
\psi_E (t_1\otimes\sigma^* u_1 )\otimes\psi_E (t_2\otimes\sigma^* u_2
)
\end{array}
\]
which induces the map
\[ \phi_E : \Lambda^2 H^0(\tC, E) \otimes \Lambda^2 H^0(\tC, \sigma^*E) \lra
S^2 H^0(\tC, \omega_{\tC})
\]
Put
\[
t := (s_1 \otimes s_2 - s_2 \otimes s_1) \otimes (\sigma^* s_1 \otimes
\sigma^* s_2 - \sigma^* s_2\otimes \sigma^* s_1)
\]
then $Z(\phi_E(t))$ is in the linear system $|\omega_{\tC^2 }|$ and is
equal to $D_E + \Delta + D_{\sigma^* E} + \Delta$. Hence $D_E +
D_{\sigma^* E}$ is in the linear system $|\omega_{\tC^2
}\otimes\cO_{\tC^2 } ( - 2 \Delta ) |$.

If $s\in\G_0$ is such that $Z(s) =\Xi_E +\Xi_{\sigma^* E}$, then, by
Lemma \ref{XiEC2}, we have $Z(\rho^* s) = D_E + D_{\sigma^* E} +
2\D'$. So $Z(\rho^* s) - 2\D' = Z(\phi_E(t)) -2\D$ and $\rho^*\cO_P (
2\Xi )\cong\omega_{\tC^2 }\otimes\cO_{\tC^2 } ( - 2 \Delta +
2\D')$. Since $Z(\rho^* s) - 2\D' = Z(\phi_E(t)) -2\D$, the section
$\rho^* s$ of $\rho^* \cO_P (2 \Xi)\otimes\cO_{\tC^2} (- 2 \Delta'
)\cong\omega_{\tC^2 }\otimes\cO_{\tC^2} (- 2 \Delta )$ is a nonzero
constant multiple of $\phi_E (t )\in I_2 (\tC )\subset H^0 (\tC^2
,\omega_{\tC^2 }\otimes\cO_{\tC^2 } ( - 2 \Delta ) )$. Since such $s$
generate $\G_0$, this proves that $\rho^*\G_0\subset I_2 (\tC )$. The
rest of part \ref{restrictC2} easily follows now.

\ref{Zs}. Let $s$ and $E$ be as above and suppose furthermore that
$E\not\cong\sigma^* E$. With the above notation, let $X_1,
X_2, X_3, X_4$ be the images of, respectively,
\[ s_1 \otimes
\sigma^* s_1, s_2 \otimes \sigma^* s_1, s_1 \otimes \sigma^* s_2, s_2
\otimes \sigma^* s_2
\]
by the map $\psi_E$. Then $\phi_E (t) = X_1X_4-X_2X_3$. By,
e.g. \cite{ACGH} page 261, the polynomial $X_1X_4-X_2X_3$ is an
equation for $q_E$. Therefore, since $\rho^* s$ is a constant nonzero
multiple of $\phi_E (t)$, we have $q(s) = q_E$. By continuity, for all
$s'$ such that $Z(s') =\Xi_M +\Xi_{\sigma^* M}$ for some $M\in\Xi'$
with $h^0 (\tC, M) = 2$, we have $q(s') = q_M$.

A line is in $q_E$ if and only if it is in a linear space of one of
the two rulings. Therefore a secant $\langle p+q\rangle$ to
$\kappa\tC$ is in $q_E$ if and only if $h^0(\tC, E\otimes\cO_{\tC} (-p-q))
> 0$ or $h^0(\tC,\sigma^* E\otimes\cO_{\tC} (-p-q)) > 0$. Hence
\[
D_E + D_{\sigma^* E} = \{ (p,q) : \langle p+q \rangle\subset q_E \} =
Z_{q_E} = Z_{q(s)}
\]
and $Z_{q(s)}$ is a reduced divisor in this case. Hence $Z_{q(s)}$ is
reduced for general $s\in\G_0$. Now the rest of part \ref{Zs} follows
by linearity.

\ref{tangcone}. When $E$ is a smooth point of $\Xi'$, the intersection
$q_E\cap\bP T_0 P = 2\bP T_E\Xi'$ is the projectivized tangent cone at $0$
to $\Xi_E +\Xi_{\sigma^* E}$ (see \cite{mumford74} pages 342-343). Hence
\ref{tangcone} also follows by linearity.

\ref{multDs}. This immediately follows from the facts that $\rho^* s\in I_2
(\tC )\subset S^2 H^0 (\tC ,\omega_{\tC })\subset H^0 (\tC^2 ,\omega_{\tC^2
})$ and that $\chi C$ is the tangent cone at $0$ to $\Si$.

\qed

\section{More divisors in $|2 \Xi |_{00}$}\label{newdiv}

For any $M \in Pic^{g-1}C$, let $\tT_{\pi^* M}$ be the translate of
$\tT$ by $\pi^* (\kappa\otimes\beta\otimes M^{-1})$, i.e.,
\[
\tT_{\pi^* M} :=\{ E\in Pic^{2g-2}\tC : h^0(\tC, E\otimes\pi^* M ) >
0\}\; .
\]
We have

\begin{proposition}\label{newdivG0}
The divisor $P.\tT_{\pi^* M}$ is in the linear system $|2 \Xi |$. It is
in $| 2 \Xi |_0$ if $h^0 (C , M)$ is positive.
\end{proposition}
{\em Proof :} We first prove that all the divisors $P.\tT_{\pi^* M}$ are
linearly equivalent as $M$ varies in $Pic^{g-1} C$. Let $\psi : JC\ra Pic^0
P$ be the morphism of abelian varieties which sends $M\otimes\kappa^{ -1
}\otimes\beta^{ -1}$ to $\cO_P (P.\tT_{\pi^* M} - P.\tT )\in Pic^0 P$. For
$C$ general, the abelian variety $JC$ is simple, hence, since its dimension
is not equal to the dimension of $Pic^0 P$, the morphism $\psi$ is trivial,
i.e., its image is $\cO_P$. By continuity, this is the case for all
$C$. Therefore, since $P.\tT = 2\Xi$, all the divisors $P.\tT_{\pi^* M}$
are linearly equivalent to $2\Xi$. The second assertion is now
immediate. \qed

\begin{proposition}\label{newdivG00}
The divisor $P.\tT_{\pi^* M}$ is an element of $| 2 \Xi |_{00}$ if $h^0(C,
M) \geq 2$.
\end{proposition}
{\em Proof :} We prove the Proposition in the case where $h^0(\tC
,\pi^* M) = 2$ and it will follow in all cases by semi-continuity.
By, e.g. \cite{ACGH} page 261, the tangent cone to $\tT_{\pi^*M}$ at
$0$ is the quadric $q_{\pi^* M} :=\cup_{\delta\in |\pi^*
M|}\langle\delta\rangle$ in $|\omega_{\tC }|^* =\bP T_0 J\tC$. Let
$\pi$ also denote the projection $\bP T_0 J\tC = |\omega_{\tC }|^*\ra
|\omega_{ C }|^* =\bP T_0 JC$ with center $|\omega_C\otimes\alpha |^*
=\bP T_0 P$. Since $q_{\pi^* M} = \pi^* (q_M ) =\pi^* (\cup_{\delta\in
| M|}\langle\delta\rangle)$, we see that $q_{\pi^* M}$ contains $\bP
T_0 P\subset\bP T_0 J\tC$ and the multiplicity of $P.\tT_{\pi^*M}$ at
$0$ is at least $3$. Since this multiplicity is even, we have
$P.\tT_{\pi^*M}\in |2\Xi |_{00}$.
\qed

\section{The base locus of $L |_{\tSi}$}\label{secproof}

In this section we prove Proposition \ref{mainprop2} and Theorem
\ref{maintheorem}. Recall that we have assumed $C$ to be
non-hyperelliptic. We have

\begin{proposition}\label{qpiD}
If $M\in Pic^{g-1} C$ is such that $h^0 ( C , M) = h^0 (\tC,\pi^* M) =
2$ and $|M|$ and $|\omega_C\otimes M^{-1 }|$ contain reduced divisors,
then
\[
\rho^* (P.\tT_{\pi^* M}) = D_M + D_{\omega_C\otimes M^{-1}} + 4\D'
\]
where $D_M$ is the reduced divisor with support
\[
\{ (p,q) : h^0 (C, M\otimes\cO_C (-\pi p -\pi q)) > 0\}\; .
\]
\end{proposition}
{\em Proof :} First note that $P.\tT_{\pi^* M}$ does not contain
$\Si$. Indeed, for general points $p$ and $q$ in $\tC$, we have $h^0 (\tC
,\pi^* M\otimes\cO_{\tC} (p+q)) = 2$ and $h^0 (\tC ,\pi^* M\otimes\cO_{\tC}
(-\sigma p -\sigma q)) = 0$. Therefore $p$ and $q$ are base points for $|
\pi^* M\otimes\cO_{\tC} (p+q) |$ and $h^0 (\tC ,\pi^* M\otimes\cO_{\tC}
(p+q -\sigma p -\sigma q)) = 0$.

Now note that $\rho^* (P.\tT_{\pi^* M}) - 4\D'$ is effective because, by
Proposition \ref{newdivG00}, the divisor $P.\tT_{\pi^* M}$ has multiplicity
at least $4$ at the origin. Then, using Lemma \ref{Zsqs}, Proposition
\ref{newdivG0} and \cite{griffithsharris78} page 287, one easily
computes that $\rho^* (P.\tT_{\pi^* M}) - 4\D'$ and $D_M +
D_{\omega_C\otimes M^{-1}}$ have the same cohomology class. It is
therefore sufficient to prove that the support of $\rho^* (P.\tT_{\pi^*
M}) - 4\D'$ contains the support of $D_M + D_{\omega_C\otimes
M^{-1}}$. By definition the support of $\rho^* (P.\tT_{\pi^* M})$ is
the set
\[
\{ (p,q) : h^0 (\tC ,\pi^* M\otimes\cO_{\tC } (p + q -\sigma p -\sigma
q) > 0\}
\]
which, by Riemann-Roch, Serre Duality and the isomorphism $\omega_{\tC
}\cong\pi^*\omega_C$, is equal to the set
\[
\{ (p,q) : h^0 (\tC ,\pi^* (\omega_C\otimes M^{-1} )\otimes\cO_{\tC }
(\sigma p +\sigma q -p-q) > 0\}\; .
\]
Therefore the support of $\rho^* (P.\tT_{\pi^* M})$ contains the union of
$\D'$ and the sets
\[
\{ (p,q) : h^0 (C, M\otimes\cO_C ( -\pi p-\pi q) ) > 0 \}
\]
and
\[
\{ (p,q) : h^0 (C,\omega_C\otimes M^{-1}\otimes\cO_C ( -\pi p-\pi q) ) > 0
\}\; .
\]
\qed
\vskip20pt

The following implies Theorem \ref{maintheorem}
\begin{proposition}\label{rhobase}
The inverse image by $\trho$ of the support of the base locus of $L
|_{\tSi}$ is contained in the set of elements $(p,q)$ of $\tC^2$ such
that $\langle\pi p + \pi q \rangle$ is contained in the intersection
of the quadrics containing the canonical curve $\kappa C$. In
particular, if $C$ is not trigonal, then the base locus of $L$ does
not intersect $\tSi$.
\end{proposition}
{\em Proof :} Let $W_{g-1}^1$ be the subvariety of $Pic^{g-1} C$
parametrizing invertible sheaves with $h^0(C,M) = 2$. By
Proposition \ref{newdivG00} the base locus of $L |_{\tSi}$ is contained in
$\tSi\cap\left(\cap_{M \in W_{g-1}^1}\epsilon^{-1}_* ( P .\tT_{\pi^* M}
)\right)$ where $\epsilon^{-1}_*( P .\tT_{\pi^* M} ) =\epsilon^* (P
.\tT_{\pi^* M}) - 4\cE$. We have $\trho^* (\epsilon^{-1}_* (P .\tT{\pi^*
M})) =\rho^* (P .\tT_{\pi^* M}) - 4\D'$. It immediately follows from
Proposition \ref{qpiD} that the inverse image $\trho^*\left(\cap_{M\in
W_{g-1}^1}\epsilon^{-1}_* ( P . \tT_{\pi^* M} )\right)$ is supported on the
set of elements $(p,q)$ of $\tC^2$ such that $\langle\pi p +\pi q\rangle$
is contained in $q_M$ for all $M\in W_{g-1}^1$. Since the quadrics of the
form $q_M$ generate $| I_2(C) |$ (see \cite{green84} and
\cite{smithvarley90}) and the base locus of $|I_2(C)|$ in the canonical
space $|\omega_C |^*$ does not contain any secants to $\kappa C$ for $C$
non-trigonal (see \cite{ACGH} page 124), the proposition follows.  \qed

\begin{remark}\label{rembase}
For $C$ trigonal, $(P,\Xi )$ is the jacobian of a curve, say $X$. Then,
using the set-theoretical equalities $V_{00} = X-X$ and $V_{inf,00} =
\kappa X$, one can show that the inverse image by $\trho$ of the support of
the base locus of $L |_{\tSi }$ is {\em equal} to the set of
$(p,q)\in\tC^2$ such that $\langle\pi p +\pi q\rangle$ is contained in the
intersection of the quadrics containing the canonical curve $\kappa C$.
\end{remark}

Let
\[
\tau_2 :\G_0\lra H^0 (\bP T_0 P,\cO_{\bP T_0 P} (2))
\]
be the map which to $s\in\G_0$ associates the quadric term of its
Taylor expansion at $0$. Then $\tau_2$ is onto because its kernel is
$\G_{00}$ which has dimension $2^p - 1 -\half p(p+1) = $ dim$(\G_0)
-h^0 (\bP T_0 P,\cO_{\bP T_0 P} (2))$. We have

\begin{lemma}
For $s\in\G_0$,
\[
s\in\G_{00}\Longleftrightarrow q(s)\supset\bP T_0 P\; .
\]
\end{lemma}
{\em Proof :} By Lemma \ref{Zsqs} part \ref{tangcone} and with the
notation there, if
$s\in\G_0\setminus\G_{ 00}$,
then $\tau_{Z(s)} = q(s)\cap\bP T_0 P$. So the projectivizations of
the two maps $\tau_2$ and $s\mapsto (\rho^* s) |_{\bP T_0 P}\in I_2
(\tC ) |_{\bP T_0 P}$ are equal. Hence, there exists $\lambda\in\bC^*$
such that, for every $s\in\G_0$, we have $\lambda\tau_2 (s) = (\rho^*
s) |_{\bP T_0 P}$. So
\[
s\in\G_{00}\Leftrightarrow\tau_2 (s) =
0\Leftrightarrow (\rho^* s) |_{\bP T_0 P} = 0\Leftrightarrow
q(s)\supset\bP T_0 P\; .
\]
\qed
\vskip20pt

Let $I_2 (\tC ,\alpha)$ be the subvector space of $I_2 (\tC)$
consisting of elements which vanish on $\bP T_0 P$. By the above
lemma and because all elements of $\G_0$ are even, the map
$\rho^*$ sends $\G_{00}$ into the subspace $I_2 (\tC ,\alpha)^+$ of
$\sigma$-invariant elements of $I_2 (\tC ,\alpha)$. We have

\begin{lemma}\label{I2inv}
The subspace $I_2 (\tC ,\alpha)^+$ is equal to $I_2 (C)\stackrel{\pi^*
}{\subset} I_2 (\tC)$.
\end{lemma}
{\em Proof :} The $\sigma$-invariant and $\sigma$-anti-invariant parts
of $H^0 (\tC ,\omega_{\tC} )$ are, respectively, $H^0 (C ,\omega_{C}
)$ and $H^0 (C ,\omega_{C}\otimes\alpha )$. Therefore in the decomposition
\[
S^2 H^0 (\tC ,\omega_{\tC} ) = S^2 H^0 (C ,\omega_{C} )\oplus H^0 (C
,\omega_{C} )\otimes H^0 (C ,\omega_{C}\otimes\alpha )\oplus S^2 H^0
(C ,\omega_{C}\otimes\alpha )
\]
the space $S^2 H^0 (C ,\omega_{C} )\oplus S^2 H^0 (C
,\omega_{C}\otimes\alpha )$ is the $\sigma$-invariant part of $S^2 H^0
(\tC ,\omega_{\tC} )$. So $S^2 H^0 (C ,\omega_{C} )$ is the subspace
of $\sigma$-invariant elements of $S^2 H^0 (\tC ,\omega_{\tC} )$ which
vanish on $\bP T_0 P$. Therefore $I_2 (\tC ,\alpha)^+$ is contained in
$S^2 H^0 (C ,\omega_{C} )$ and $I_2 (\tC ,\alpha)^+$ is the subspace
of elements of $S^2 H^0 (C ,\omega_{C} )$ which vanish on
$\kappa\tC$. This is precisely $I_2(C)$. \qed

\begin{corollary}
The dimension of $\G_{00}'$ is at least $2^p - 2 - p^2 +p$ and at most
$2^p - 1 -\half p(p+1) - 3 = 2^p - 4 -\half p(p+1)$. The codimension
of $\cQ_{00}'$ in $\cQ_{00}$ is at least $2$.
\end{corollary}
{\em Proof :} By Lemma \ref{I2inv}, we have $\G_{00}' = Ker (\rho^*
:\G_{00}\ra I_2 (C)\stackrel{\pi^* }{\subset} I_2 (\tC))$. Hence the
dimension of $\G_{00}'$ is at least dim$(\G_{00}) -$ dim$(I_2 (C)) =
2^p - 1 -\half p(p+1) -\half (g-2)(g-3) = 2^p - 1 -\half p(p+1) -\half
(p-1)(p-2) = 2^p - 2 - p^2 +p $.

To prove the upper bound for the dimension of $\G_{00}'$, we prove
that the dimension of $\rho^*\G_{00}$ is at least $3$. For $C$
non-trigonal, this is an immediate consequence of the result of
Proposition \ref{rhobase} which says that $\trho^* L$ has no base
points and the fact that $\trho^* L$ has positive self-intersection
number: this number is easily computed using Lemma \ref{XiEC2}
and the calculation of the cohomology class of $D_E$ on page 287 of
\cite{griffithsharris78}.

Now suppose $C$ trigonal. The linear system $\trho^* L$ contains the
divisors $\trho^* (P.\tT_{\pi^* M}) - 4\D'$ for $M\in
W^1_{g-1}$. Since $g\geq 5$, the curve $C$ has a unique linear system
of degree $3$ and dimension $1$ and we denote the associated
invertible sheaf of degree $3$ by $M_0$. Choose three general
effective divisors $N$, $N'$ and $N''$ of degree $g-4$ on $C$ and put
$M = M_0\otimes\cO_C( N)$, $M' = M_0\otimes\cO_C( N')$ and $M'' =
M_0\otimes\cO_C( N'')$. Let $s$, $s'$ and $s''$ be sections whose
divisors of zeros are respectively $\trho^* (P.\tT_{\pi^* M}) - 4\D'$,
$\trho^* (P.\tT_{\pi^* M'}) - 4\D'$ and $\trho^* (P.\tT_{\pi^* M''}) -
4\D'$. If a non-trivial linear combination of $s, s', s''$ is zero,
then the coefficients of at least two of $s, s', s''$, say $s$ and $s'$,
are non-zero. Then the restrictions of $s$ and $s'$ to $\trho^*
(P.\tT_{\pi^* M''}) - 4\D'$ have the same divisors. By Proposition
\ref{qpiD}, we have $\trho^* (P.\tT_{\pi^* M}) - 4\D' = D_{M}
+D_{\omega_C\otimes {M}^{-1}}$. Hence $\trho^* (P.\tT_{\pi^* M}) -
4\D' = D_{M_0} + p_1^* (\pi^* N) +p_2^* (\pi^* N) +D_{\omega_C\otimes
{M}^{-1}}$. Similarly, $\trho^* (P.\tT_{\pi^* M'}) - 4\D' = D_{M_0} +
p_1^* (\pi^* N') +p_2^* (\pi^* N') +D_{\omega_C\otimes {M'}^{-1}}$ and
$\trho^* (P.\tT_{\pi^* M''}) - 4\D' = D_{M_0} + p_1^* (\pi^* N'')
+p_2^* (\pi^* N'') +D_{\omega_C\otimes {M''}^{-1}}$. Let $p$ be a
point of $\tC$ such that $\pi p$ is a point of $N''$. Then $F_p :=
p_1^* (p)\cong\tC$ is contained in $\trho^* (P.\tT_{\pi^* M''}) -
4\D'$. Hence $(\trho^* (P.\tT_{\pi^* M}) - 4\D' )|_{F_p} = (\trho^*
(P.\tT_{\pi^* M'}) - 4\D' )|_{F_p}$. Equivalently, $(p_2^* (\pi^* N) +
D_{\omega_C\otimes {M}^{-1}} )|_{F_p} = (p_2^* (\pi^* N') +
D_{\omega_C\otimes {M'}^{-1}} )|_{F_p}$.  Since $N$, $N'$, $N''$ are
general, there are unique divisors, say $G$ and $G'$ in, respectively,
$|\omega_C\otimes {M}^{-1} |$ and $|\omega_C\otimes {M'}^{-1} |$
containing $\pi p$. Then, after identifying $F_p$ with $\tC$, the
previous equality of divisors becomes $\pi^* N + \pi^* G -p-\sigma p
=\pi^* N' + \pi^* G' -p-\sigma p$. Since $N'$ is general and $G$
depends only on $N$ and $N''$, this is not possible. Hence there are
no non-trivial linear relations between $s, s', s''$ and $\trho^* L$
is at least a net.

In the non-trigonal case, the codimension of $\cQ_{00}'$ in $\cQ_{00}$
is at least $2$ because the restriction of $\trho^* L$ to $\D'$ has
positive degree and no base points, therefore it is at least a pencil.
In the trigonal case, use divisors as above to show that $(\trho^*
L)|_{\D'}$ contains at least two distinct divisors.
\qed

\begin{remark}
In case $p=4$, the above gives a second proof of the result of
\cite{I3} (page 148) saying that the dimension of $\G_{00}'$ is $2$.
\end{remark}


\providecommand{\bysame}{\leavevmode\hbox to3em{\hrulefill}\thinspace}

\end{document}